\documentclass[11pt]{article}
\usepackage{amsmath}
\usepackage{amssymb}
\usepackage{hyperref}

\usepackage{graphicx}
\usepackage{comment}
\usepackage{xcolor}

\usepackage{amsthm}
\theoremstyle{definition}
\newtheorem*{Definition1}{Definition   1}
\newtheorem*{Theorem1}{Theorem   1}
\newtheorem*{Theorem2}{Theorem   2}

\newtheorem*{pthm2}{Proof of  Theorem  2}
\newtheorem*{Q}{Question    1}
\newtheorem*{QQ}{Question   2}
\newtheorem*{QQQ}{Question  3}

\newtheorem*{pro1}{Proposition    2}
\newtheorem*{Lemma}{Lemma}

\begin{document}
\title{A Note on Application of  Singular rescaling  to Global Bifurcation of Lienard Equations}
\author{Ali Taghavi \\ Qom University of Technology}
\maketitle
\begin{abstract}
We apply the  singular rescaling $x:=\frac{x}{\epsilon},\;y:=\frac{y}{\epsilon}$ to give an  alternative proof  for a result in local bifurcation of  Lienard equations proved in A. Lins, W. de Melo and C. C. Pugh, \emph{On Lienard's Equations},
\textit{lecture notes in Mathematics, 597. Springer verlag(1977)}. We also give a generalization of this  local result to a global one.
\end{abstract}
\section{Introduction}
 Hilbert  16th problem asks for a uniform upper bound $H(n)$ for the number of  limit cycles of a polynomial vector  field \begin{equation*} \begin{cases}  x'=P_n(x,y)\\y'=Q_n(x,y) \end{cases} \end{equation*}
where $P_n,Q_n$  are polynomials in $\mathbb{R}[x,y]$. The  problem  among a list of 23 problems has been posed by David Hilbert  in 1900  congress of  mathematician in Paris. The problem list is published  in \cite{Hilbert}. The  Hilbert 16th problem is still an open problem  in its  full generality. A weaker version of this problem is  to  study the number of  limit  cycles which appear by  small perturbation of  a  given  vector  field $X_0$. More  precisely assume that we  have a finite parameter family $X_\lambda$ of  vector  fields  where $\lambda$ is a finite dimensional parameter. We are interested in the  maximum number of  limit  cycles of $X_\lambda$  for  $\lambda$  sufficiently close to $\lambda=0$. Recall that a limit cycle is an isolated periodic orbit for a 2 dimensional vector field.
To  have a better study of dynamical behaviour of  a polynomial system on $\mathbb{R}^2$ one  transfers  a given  planar  system  to the Poincare sphere $S^2$. The process is called Poincare  compactification. The study of bifurcation of  limit cycles under small perturbation of parameter in a   finite parameter family of vector fields  generates the  concept of cylicity  of limit  periodic  sets mentioned in \cite{Roussarie}. The  definition is as  follows:

\begin{Definition1}
Let $X_\lambda$ be  a family of vector fields on $\mathbb{R}^2$  where  $\lambda$ varies in a  subset $D$ of $\mathbb{R}^k$. The family is called a $k$ parameter family of vector  fields and $D$ is called the parameter space. Let $\lambda_n\in D$ converges to $\lambda_0$ and  $X_{\lambda_n}$ possess a limit cycle $\gamma_n$. Assume that $\gamma_n$ converges to a compact subset $\gamma$ of $\mathbb{R}^2$ with respect to the  Hausdorff metric   .  We call $\gamma$ a limit periodic set. The  cyclicity of $\gamma$ is defined as the  maximum number of limit cycles of $X_\lambda$ which are Haussdorf close to $\gamma$ for $\lambda$   sufficiently close to $\lambda_0$.
\end{Definition1}

It is well known that a limit periodic set is invariant under flow of $X_0$. When $X_0$ has a finite number of  singularities then Poincare Bendixon theorem characterise all limit periodic sets: every limit periodic set is either a periodic orbit or a separatrix polycycle consist of a finite (possibly a single) number of  singularities and (possibly an empty set of) regular oriented orbits joining the singularities. A  separatrix cycle is  called a graph. Hence a  single  singularity is  counted as a graph thought there is  no any regular orbit joining it to other (or itself) singularity.

An important class of  polynomial vector fields is the Lienard equation:
\begin{equation}\label{Main}
\begin{cases}
x'=y-F(x)\\
y'=-x
\end{cases}
\end{equation}
where $F(x)$ is a polynomial function  with $deg(F)=2n+1$ or $deg(F)=2n+2$.  It is  conjectured in \cite{3} that the  maximum number of  limit cycles of the system \eqref{Main} is $n$. The  conjecture is answeredno by negative in \cite{Roussarie2}The reason for belief on the  conjecture was based on two types of bifurcation results. The first one was that the order of a weak focus  in \eqref{Main} is at most $n$ so a weak focus of this family can generate at most $n$  small  limit cycles around the focus singularity. Recall that a weak focus is  a  singularity of  a  vector  field whose linear part is in the form \begin{equation*} \lambda y \partial_x-\lambda x \partial_y \end{equation*} Locally around a  weak focus singularity one defines a Poincare return map $P$  with  an isolated fixed point at the origin. The order of  weak focus is the order of origin as an isolated root  for the function $P-Id$.  With  a  straightforward computation one may compute the order of  weak focus in Lienard family. The method of such computations date backs to the time of publication of seminal paper  \cite{Bautin}. See \cite{4} for  a  survey of  such kind of  weak focus computations. The precise focal result for the Lienard family is the following:
\begin{pro1}
The origin is a  weak  focus of order $n$  for $\begin{cases} x'=y-x^{2n+1}\\ y'=-x \end{cases}$
\end{pro1}

 The other reason for  belief on the conjecture is the following  result on the number of zeros of  Abelian integral proved in \cite{3}:

\begin{Theorem1}
Let $F(x)=a_{2n+1}x^{2n+1}+a_{2n}x^{2n}+\ldots+a_1 x$ be a polynomial of degree $2n+1$. Then the linear center   \begin{equation*} \begin{cases} x'=y\\y'=-x \end{cases} \end{equation*} as unperturbed system of one parameter family \begin{equation}\label{LIENARD} \begin{cases}x'=y-\epsilon F(x)\\y'=-x \end{cases} \end{equation}  has at most n circle $\gamma_1, \gamma_2,\ldots,\gamma_n$ as  limit periodic sets.
\end{Theorem1}

The proof of the theorem is based on averaging method and counting the number of  zeroes of Abelian integral $I(h)=\oint_{\gamma_h} F(x)dy$  defined for $h\in (0, \infty)$ where $\gamma_h=H^{-1}(h)$ with $H(x,y)=x^2+y^2$. The Hamiltonian $H$ is a first integral for the linear center \begin{equation*} \begin{cases} x'=y\\y'=-x \end{cases} \end{equation*} Note that the theorem  above does not imply  that the  system \eqref{LIENARD} with sufficiently small $\epsilon$ has at most $n$ limit cycle in whole plane. But it  states  that there are at most n circles which may generate limit cycles. In  fact the theorem does not have any control on the number of limit cycles which escape at infinity or escape at origin for  sufficiently small parameter $\epsilon$. The theorem merely control the limit cycles in the  compact subsets of the punctured plane $\mathbb{R}^2\setminus \{0\}.$  In this paper we  give  the following generalization of the above theorem as  follows:

\begin{Theorem2}
For  sufficiently small $\epsilon $ the system \eqref{LIENARD} has at most $n$ limit cycles in the whole plane.
\end{Theorem2}

Theorem 2 is stronger than Theorem 1. This gives us in particular an alternative  proof for the result in \cite{3} when the degree of  polynomial $F(x)$ is odd.  As we shall see in the next section the proof of Theorem 2 is not based  on application of  averaging method and  Abelian integrals.

\section{Preliminaries and Proofs}
In this paper a planar  vector  field $\begin{cases} x'=P(x,y)\\y'=Q(x,y)  \end{cases}$ is some times denoted by  $P\partial_x+Q\partial_y$. A  limit cycle  for a  planar  vector  field is an isolated periodic orbit. A limit cycle $\gamma$ is  called a stable   limit  cycle if  the solutions starting nearby point to $\gamma$ tend towards $\gamma$ as parameter $t$ goes to infinity. It is  unstable  if the  solutions starting nearby points tends to $\gamma$ as  $t\to -\infty$. It is called  a semi stable  limit cycle if it is  stable  from interior(or  exterior) and  unstable from exterior(interior). It is well known that a limit cycle which is not a semi stable limit cycle persist under  small perturbation of the  vector  field. Moreover a  semi stable limit cycle would be twisted (or disappear)  with an appropriate (or unappropriate)perturbation of the  underling vector field, see \cite{Perko}.\\
We shall  introduce the Poincare compactification method to analysis the  dynamics of a polynomial vector field.
Let we have a Polynomial vector field  $X=(P_n(x,y),Q_n(x,y)$ of degree $n$ on the plane. We transfer  $X$ to the unit sphere $S^2=\{(x,y,z)\in \mathbb{R}^3 \mid  x^2+y^2+z^2=1\} $. The upper  and  lower  hemispheres  are denoted by $S^2_+$ and $S^2_-$ respectively. We consider two diffeomorphisms $\phi_{\pm}:\mathbb{R}^2 \to S^2_{\pm}$  with
$$\phi_{\pm}(x,y)=\left (\frac{x}{\sqrt{1+x^2+y^2}},\frac{y}{\sqrt{1+x^2+y^2}},\frac{\pm1}{\sqrt{1+x^2+y^2}}\right)$$
The push forward vector fields $\phi_{\pm}*(X)$ on upper  and lower hemi sphere can not be  extended to the equator $z=0$ of $S^2$ since they are unbounded vector fields near the equator. But after rescalling by  term $z^{n-1}$ we can extend $z^{n-1}\phi_{\pm}*(X)$ to an analytic vector field $\tilde{X}$ on $S^2$. The vector field $\tilde{X}$ is called the Poincare compactification of polynomial vector  field $X$. In this way the infinity of the plane is  compactified to the  equator $z=0$ on the Poincare sphere $S^2$. So those orbits in the plane who scape to infinity, namely unbounded orbits,  would accumulate  at equator so we can analysis the behaviour of such orbits by study  their accumulation points on the equator. The phase portrait of the Poincare compactification of the polynomial  Lienard equation \begin{equation}\label{LIENARDD} \begin{cases} x'=y-F(x)\\y'=-x\end{cases} \end{equation} is completely described in \cite{3}.   \\

Let origin $(0,0)\in \mathbb{R}^2$ be  a  singular point of  an analytic vector field on a neighborhood of the origin whose linear part is in the form \begin{equation} \begin{cases}x'=y\\y'=-x \end{cases}\end{equation} If the origin is surrounded by a  band of closed orbit then it is called a center. Otherwise it is called a  weak focus. One can define a Poincare return map $P$ defined on a local section parameterized by $s\in (-\delta, \delta)$  around a weak  focus. The order of a weak focus is the order of 0 as a root of the function $P-Id$. It is  well known that the ciclicity of  a  weak  focus of order $k$ in a  finite parameter family of vector  fields is at most $k$. So with  small perturbation of a weak focus we get at most $k$ small limit cycles locally around the singularity, see \cite{9} and \cite{Perko}.\\

We equip the space $\chi^\infty(S^2)$ of  smooth vector fields on $S^2$ with the  compact open topology or  weak or  strong topology. The latter two topologies are the same since the  underling space $S^2$ is a compact space. For  definitions of these  topologies  \cite{Hirsch}.

The  following lemma is  a consequence of definition of Poincare  compactification.

\begin{Lemma}
Let $X_\lambda$ be a  $k$ parameter family of  polynomial vector fields of degree $n-1$.Assume that $\lambda_n $ tends to $0\in \mathbb{R}^k$. Moreover assume that $Y$ is a polynomial vector field of degree $n$. Then the  Poincare compactification of $X_{\lambda_n}+Y$ tends to the  compactification of $Y$ in the compact open topology or weak-strong topology.
\end{Lemma}

The processes we mentioned above and the description of  phase portrait of \eqref{Main} described in \cite{3}  leads us to the following  proposition:

\begin{pro1}
Assume that $g(x,\epsilon)$ is a smooth function which is a polynomial function of degree $2n$ in variable $x$. Then the  only limit periodic set in $S^2$ of $\begin{cases} x'=y-ax^{2n+1} \\ y'=-x \end{cases}$ in the family $\begin{cases}\label{Perturbation} x'=y-ax^{2n+1}+g(x,\epsilon) \\ y'=-x \end{cases}$ is $(0,0,\pm 1)$

\end{pro1}



Now we are ready to give a proof  of Theorem 2 which generalizes Theorem 1 in \cite{3}:
\begin{pthm2}
We put \begin{equation}\label{singular} \epsilon:=\delta^{2n}, x:=\frac{x}{\delta}, y:=\frac{y}{\delta}\end{equation} After this change of  coordinate we get the following equation:
\begin{equation}\label{EQEQ}  \begin{cases} x'=y- a_{2n+1}x^{2n+1}-\left(\delta a_{2n}x^{2n}+\delta^2 a_{2n-1}x^{2n-1}+\ldots+\delta^{2n}a_1 x\right )\\y'=-x \end{cases}\end{equation}
The  unperturbed  system  with $\delta=0$ is \begin{equation}\label{12345}\begin{cases}x'=y-a_{2n+1}x^{2n+1}\\y'=-x  \end{cases}\end{equation} So the system \eqref{EQEQ}  is a perturbation of a vector field of degree $2n+1$ with small terms of degree 2n. So this is a perturbation in $\chi^{\infty} (S^2)$  with respect to compact open topology So limit cycles of \eqref{EQEQ} merge in a limit periodic set of \eqref{12345}. By proposition 2 the only limit  periodic set is a weak focus of order n hence it can generate at most $n$ limit cycles. This implies that for  sufficiently small $\delta$ \eqref{EQEQ} has at most $n$ limit cycles in the whole plane. This  completes the  proof of the  theorem.
\end{pthm2}

\section{Conclusion and  some questions for further research} In this paper we considered the  equation  $$ \begin{cases}x'=y-\epsilon F(x)\\y'=-x \end{cases}$$ where $F(x)$ is a  polynomial of degree $2n+1$. We used the singular rescalling  $$x:=\frac{x}{\delta},\quad y:=\frac{y}{\delta}, \quad\epsilon=\delta^{2n}$$ to prove that the  above system  has at most $n$ limit cycles in the whole plane. This  is  a stronger result compared to a local version proved in \cite{3}. We observed that this  perturbational result can be proved in a stronger version independent of theory of  Abelian integrals  because the  focus problem already contained all intrinsic information.  Our global results does not work if the polynomial $F(x)$ would be an even degree polynomial. When $F(x)=x^{2n}+\text{Lower degree terms}$ then  after the  above  singular rescalling we  get a  perturbation of the system $\begin{cases} x'=y-x^{2n}\\ y'=-x\end{cases}$ which has a  center at origin and  a  homoclinic  loop based  at points $(0,1,0)$ at infinity. To our knowledge the  cyclicity of this  homoclinic loop and  also the  Abelian integrals associated to the perturbation of this system with center is not investigated yet. Interestingly enough the unperturbed system is  not  a  Hamiltonian system and  does not  admit a first integral with an explicit formula in quadratures.\\
In the paper we applied the  singular rescaling to a  particular  algebraic    vector  field namely the  Lienard vector  field. One can examine this  rescaling to arbitrary polynomial vector  field to obtain some  perturbation result and possibly an upper bound for the number of  limit cycles of  such systems. One observe that the  processers  does not work in general case. But as a byproduct of this processes we produce two  questions apparently unrealted to dynamics  but  are in the context of  Lie  algebras and  $C^*$ algebras. \\

 Let $H(n)$ be  the uniform upper bound for the number of limit cycles of a polynomial vector field of degree $n$. So we would  have  a  polynomial vector field $X=P_n\partial_x+Q_n\partial_y$ of degree $n$  which possess at least  $H(n)$ limit cycles. From methods of  bifurcation theory we may assume that non of these limit cycles are semi stable limit cycle. So we conclude that the system  $X_\lambda=(P_n+\lambda x(x^2+y^2)^{n})\partial_x+(Q_n+\lambda y(x^2+y^2)^n)\partial_y$ has at least $H(n)$ limit cycles  when $\lambda$ is  sufficiently small. After a linear rescalling we may assume that $\lambda=1$ so $P_n,Q_n$  would be replaced with new polynomial functions but we denote  them with $P_n, Q_n$ again. So we  have  a  polynomial vector field in the form
  \begin{equation}\label{37} \begin{cases} x'=P_n(x,y)+x(x^2+y^2)^n\\ y'=Q_n(x,y)+y(x^2+y^2)^n\end{cases}\end{equation} with at least $H(n)$ limit cycles. We apply the  singular rescalling $x:=\frac{x}{\epsilon},\quad y:=\frac{y}{\epsilon}$ to \eqref{37} then multiply with term $\epsilon^{2n}$.  We  get a perturbation \begin{equation}\label{25} \begin{aligned} X_\epsilon&=\left(\epsilon^{2n+1}P_n(\frac{x}{\epsilon},\frac{y}{\epsilon})+x(x^2+y^2)^n\right)\partial_x+\\
  &\left(\epsilon^{2n+1}Q_n(\frac{x}{\epsilon},\frac{y}{\epsilon})+y(x^2+y^2)^n\right) \partial_y  \end{aligned}  \end{equation} of the  homogeneous vector field \begin{equation}\label{24} X_0=\left(x(x^2+y^2)^n\right)\partial_x+\left(y(x^2+y^2)^n\right) \partial_y   \end{equation} which perturbational terms are polynomials of degree at most $n$. On the other hand all limit cycles of the  systems \eqref{37} would shrink to the origin via the inverse transformation $x:=\epsilon x,\quad y:=\epsilon y$ of the above singular rescalling. So $X_\epsilon$ has at least $H(n)$  small limit cycles shrinking to the origin. However the origin as limit  periodic set of  \eqref{24}  is  a quite degenerate  singularity but it would be interesting to compute its cyclicity in some particular case. Note that $[X_0, Y]=0$ where $Y=y\partial_x-x\partial_y$ is linear center and  $X_0$ is the  system \eqref{24}. Recall from previous section that for  every two  vector fields  $X,Y$ with $[X,Y]=0$ every limit cycle of  $X$ is invariant under flow of $Y$. So in order to study the number of limit cycles of \eqref{37} and its equivalent form \eqref{24} it is natural to ask: can one  generate a  one  parameter family  $Y_\epsilon$ as a perturbation of the linear center $Y$  with $[X_{\epsilon},Y_{\epsilon}]=0$? As we said   the latter equation implies that  every limit  cycle of  $X_{\epsilon}$ is  invariant under $Y_{\epsilon}$. So the problem converts to study of the  cyclicity of origin under perturbation of a linear center. This  cylicity  denoted by $C(n)$ is finite by the concept of  Bautin ideal mentioned in \cite{Bautin} and
\cite{Roussarie}. So if the answer to the latter question is affirmative then the  process above  implies that $H(n)$ is less than $C(n)$ hence is  finite. In reality the  answer to this  question is  obviously negative   because all singularities of \eqref{37} would shrink to the origin. This would implies that a perturbation of linear center possess more than one singularity near the origin which is impossible  since every perturbation of  linear center  with higher degree polynomials has an isolated  singularity at origin. However this  question has an obvious negative answer but according to its  lie  algebraic nature we may generate some questions in the thoory of  Lie  algebras and $C^*$ algebras. Note that in the mentioned process we  had  two vector  field  $X,Y$   with  $[X,Y]=0$. For a  perturbation $X_\epsilon$ of $X$ we  wish  to  find a perturbation $Y_\epsilon$ of  $Y$  with $[X_\epsilon, Y_\epsilon]=0$. We  simulate this  situation in an  arbitrary  finite  dimensional Lie  algebra as follows:\\

Let  we have  a finite dimensional  Lie  algebra $L$  with two  non zero commuting  elements $[a,b]=0$. Are there neighborhoods $U\ni a, V\ni b$  and  smooth  function $\phi:U\to V$ with $\phi(a)=b, [x,\phi(x)]=0$? For  Abelian  Lie  algebra or cross product Lie  algebra structure on $\mathbb{R}^3$ the answer is affirmative. For  the  Matrix  algebra the  answer  is negative. In the  matrix  algebra  case there is  no  even  a  continuous  $\phi$ with the above  mentioned property. In fact  when we try to give  a continuous  $\phi$ via  Micheal Selection theorem  we realize that the  set  valued function $\phi:M_n\to  M_n$ with $\phi(x)=C(x)$ is  not  hemi lower  continuous.

This is  a  motivation to posses the  following  questions:\\
\begin{Q}
What is  a  complete classification of   all finite dimensional Lie  algebras for  which the  mapping $\phi:L \to L$ with $\phi(x)=C(x)$  is a hemi lower  continuous?

\end{Q}

The  above  question can not  be stated  in an infinite dimensional  Lie  algebra without  extra  continuous  structure. So  we need to  have some  topological structure. In the complex  setting one possible  choice is  a $C^*$ algebra. Recall that a $C^*$ algebra is a a complex  algebra  with an involution $*$ and  a norm satisfying $|xx^*|=|x|^2,  |xy|\leq |x||y|$  for  all $x,y$ in the algebra.  So we pose the following two questions:

\begin{QQ}
Is a  unital  algebra  with  hemi lower  continuous map  $\phi_A$ necessarily necessarily a  commutative  algebra?
\end{QQ}

\begin{QQQ}
What is  an  example of  a  simple  (non unital) algebra for  which $\phi_A$ is  hemi lower continuous?

 \end{QQQ}

\noindent

\begin{thebibliography}{30}

\bibitem{Bautin}  N.N. Bautin, \emph{On the number of limit cycles appearing with variation of the coefficients from an equilibrium state of the type of a focus or a center}, Matematicheskii Sbornik. Novaya Seriya, 1952, Volume 30(72), Number 1, Pages 181–196

\bibitem{Hilbert} David  Hilbert, \emph{ Mathematical Problems}. Bulletin of the American Mathematical Society. 8 (10): 437–479.(1902)


\bibitem{Hirsch}  M.W. Hirsch, \emph{Differential Topology} Springer, New York, Heigelberg, Berlin, 1976.

\bibitem{Roussarie2}F.Dumortier,D. Panazzolo, R. Roussarie,\emph{More Limit Cycles than Expected in Liénard Equations},Proc.Amer. math. Soc (135) 6, 2007

\bibitem{Roussarie}R. Roussarie\emph{Bifurcations of Planar Vector Fields and Hilbert's Sixteenth Problem},Birkhauser, 1998
\bibitem{3}
A. Lins, W. de Melo and C. C. Pugh,\emph{ On Lienard's Equations},
\textit{lecture notes in Mathematics, 597. springer verlag(1977)}
\bibitem{Perko}
K. M. Perko, Differential Equations and Dynamical systems,
springer, (2001)
\bibitem{9}
Zhang Zhi-Fen, et al, qualitative theory of Differential Equation,
\textit{Amer. Math. Soc. Providence(1991)}
W. A. Copple, A survey of quadratic systems, \textit{J.
Differential Equation 2(1966) 293-304.}

\end{thebibliography}
\end{document}